\documentclass[12pt]{amsart}

\usepackage{amsmath,amsthm, amssymb}
\usepackage[T1]{fontenc}
\usepackage{latexsym}
\usepackage{subfigure, epsfig, epsf}
\usepackage{color,fancyhdr,lastpage,graphicx}
\usepackage{xspace}
\usepackage{wasysym}
\usepackage{setspace}
\usepackage{bbm}
\usepackage{stmaryrd}

\usepackage{pstricks}
\usepackage{pst-node}
\usepackage{pst-tree}

\newcommand{\RR}{\mathbb{R}}        
\newcommand{\MMM}{\mathbb{M}}

 \newcommand{\mcC}{\mathcal{C}}

    \newcommand{\bfX}{{\bf X}}

\newcommand{\bbX}{\mathbb{X}}

\newcommand{\EE}{\mathbb{E}}

\newcommand{\ep}{\epsilon}
\renewcommand{\leq}{\leqslant}
\renewcommand{\geq}{\geqslant}

\newcommand{\wrt}{with respect to }

\newcommand{\ssk}{\smallskip}

\newtheorem{thm}{\hspace{-0.15cm}  {\sc Theorem} }
   \newtheorem*{definition}{ \hspace{-0.3cm} {\sc Definition}}

\numberwithin{equation}{section} 

\newenvironment{DemThmRegularity}{%
    \begin{list}{\hspace{0.5cm}{\sc Proof of theorem \ref{ThmRegularity} --}}{%
        \setlength{\topsep}{0pt}%
        \setlength{\leftmargin}{0pt}%
        \setlength{\rightmargin}{0pt}%
        \setlength{\listparindent}{0pt}%
        \setlength{\itemindent}{0pt}%
        \setlength{\parsep}{0pt}%
        \addtolength{\leftmargin}{20pt}%
        \addtolength{\rightmargin}{0pt}%
    } \item }{\hfill{\space $\rhd$}\end{list}\smallskip}

\title[Regularity of the It\^ o-Lyons map]{Regularity of the It\^ o-Lyons map}
\date{\today}
\author{I. Bailleul}
\address{IRMAR, 263 Avenue du General Leclerc, 35042 RENNES, France}
\email{ismael.bailleul@univ-rennes1.fr}

\begin{document}

\maketitle

\begin{abstract}
We show in this note that the It\^o-Lyons solution map associated to a rough differential equation is Fr\'echet differentiable when understood as a map between some \textit{Banach spaces} of controlled paths. This regularity result provides an elementary approach to Taylor-like expansions of Inahama-Kawabi type for solutions of rough differential equations depending on a small parameter, and makes the construction of some natural dynamics on the path space over any compact Riemannian manifold straightforward, giving back Driver's flow as a particular case.
\end{abstract}

\section{Introduction}
\label{SectionIntroduction}

It is probably fair to say that, from a probabilist's point of view, the main success of the theory of rough paths, as developed originally by T. Lyons \cite{Lyons97}, was to provide a framework in which a notion of integral can be defined as a \textit{continuous} function of both its integrand and its integrator, while extending the It\^o-Stratonovich integral when both theories apply. This continuity result is in striking contrast with the fact that the stochastic integration map defines only a measurable function of its integrator, with no hope for a better dependence on it as a rule. The continuity of the rough integration map provides a very clean way of understanding differential equations driven by rough signals, and their approximation theory. However, the highly nonlinear setting of rough paths and its purely metric topology prevent the use of the classical Banach space calculus in this setting; so, as a consequence, one cannot hope for a better statement than the following one for instance: Under some conditions to be made precise, the solution to a rough differential equation is a locally Lipschitz continuous function of the driving rough path. 

\ssk

It is fortunate that Gubinelli developped in \cite{Gubinelli} an alternative approach to rough differential equations based on the \textit{Banach space} setting of paths controlled by a fixed rough path, and which somehow allows to linearize many considerations. (We refer the reader to the book \cite{FH13} of Friz and Hairer for an excellent short account of rough path theory from Gubinelli's point of view.) We show in this note that the It\^o-Lyons solution map that associates to some Banach space-valued controlled path $y_\bullet$ the solution to a rough differential equation driven by $y_\bullet$ is actually a Fr\'echet regular map of both the controlled path and the vector fields in the equation. This regularity result, theorem \ref{ThmRegularity} in section \ref{SectionRegularity}, provides a straightforward approach to investigating the dependence of the solution to a parameter-dependent rough differential equation as a function of this parameter (section \ref{SubSectionTaylorExpansion}), as in the works of Inahama and Kawabi \cite{InahamaKawabi}, and to constructing some dynamics on some path spaces in a geometrical setting (section \ref{SubSectionDynamicsPathSpace}), such as Driver's flow.

\medskip

We have chosen to present our results in the setting of paths controlled by a $p$-rough path, with $2\leq p<3$. This makes the use of controlled paths friendly, avoiding the use of branched rough path whose algebraic structure may seem complicated to some readers, see \cite{GubinelliBranched}, while it will be clear for those acquainted with branched rough paths that nothing has to be changed to deal with the general case. 

\bigskip

\noindent {\bf Notations.} A few notations will be used throughout the note, which we gather here.

\begin{itemize}
   \item We denote by E, H, U and V some Banach spaces, and by $\textrm{L}(\textrm{U},\textrm{V})$ the set of continuous linear maps from U to V, endowed with the operator norm.  \vspace{0.1cm}
   
   \item For a function $f : \textrm{U}\rightarrow\textrm{H}$, of class $C^k$, we denote by $f^{(k)}$ its $k^\textrm{th}$-derivative and write
   $$
   \big\|f^{(k)}\big\| := \underset{u\in\textrm{U},\,|u|=1}{\sup}\;\big|f^{(k)}(u,\dots,u)\big|_\textrm{H}   
   $$
   for the operator norm of $f^{(k)}$. (As $f^{(k)}$ is symmetric for $k\geq 2$, this definition of the operator norm is equivalent to the usual one.) \vspace{0.1cm}
   
   \item Given a positive real number $\gamma$, we denote by $\textrm{Lip}_\gamma(\textrm{U},\textrm{H})$, resp. $\textrm{Lip}_\gamma(\textrm{U})$, the set of H-valued, resp. real-valued, maps on U that are $\gamma$-Lipschitz in the sense of Stein; these maps are in particular bounded. The $\gamma$-Lipschitz norm of $f\in \textrm{Lip}_\gamma(\textrm{U},\textrm{H})$ is defined by the formula
   $$
   \|f\|_\gamma := \sum_{k=0}^{[\gamma]}\big\|f^{(k)}\big\| + \underset{x\neq y}{\sup}\; \frac{\big\|f^{(k)}(x) - f^{(k)}(y)\big\|}{\|x-y\|^{\gamma-[\gamma]}}.
   $$
   If $0<\gamma<1$, the set of $\gamma$-Lipschitz functions is simply the set of $\gamma$-H\"older functions.   \vspace{0.1cm}
   
   \item Given $k\geq 3$, a continuous linear map F from U to the space of  $\textrm{Lip}_k$ vector fields on V is called a $\textrm{Lip}_k$ \textbf{V-valued $1$-form on U}. The natural operator norm on the space $\textrm{L}\big(\textrm{U},\textrm{Lip}_k(\textrm{V},\textrm{V})\big)$ of all such maps turns this space into a Banach space.  \vspace{0.1cm}  

   \item We shall denote by $\textrm{E}^{\otimes 2}$ the completion of the algebraic tensor product $\textrm{E}\otimes_a\textrm{E}$, \wrt a tensor norm for which $\textrm{L}\big(\textrm{E},\textrm{L}(\textrm{E},\textrm{H})\big)$ is continuously embedded in $\textrm{L}\big(\textrm{E}^{\otimes 2},\textrm{H}\big)$. The injective tensor norm has this property for instance. Similar assumptions are made on the tensor products $\textrm{V}^{\otimes 2}$ and $\textrm{H}^{\otimes 2}$ that appear below.
\end{itemize}

\medskip

\section{Regularity of the It\^o-Lyons solution map}
\label{SectionRegularity}

We prove in this section that the It\^o-Lyons solution map to a rough differential equation is Fr\'echet regular, when properly defined on some space of controlled paths; this is theorem \ref{ThmRegularity} below. We first recall for the reader's convenience the basics about controlled paths.

\bigskip

The starting point of Gubinelli's approach to rough paths theory is the fact that an E-valued path $x_\bullet$ over some finite interval $[0,T]$, has increments $x_{ts} := x_t-x_s$, for $s\leq t$, that enjoy the additivity property
\begin{equation}
\label{EqAdditivity}
x_{ts}-\big(x_{tu}+x_{us}\big) = 0
\end{equation}
for all $0\leq s\leq u\leq t\leq T$, and that an E-valued $2$-index map $\mu := (\mu_{ts})_{0\leq s\leq t\leq T}$ for which the identity \eqref{EqAdditivity} holds gives the increments of a path $x_\bullet$ uniquely determined, up to its initial value. Better, if $\mu$ is almost additive, in the sense that one has 
$$
\big|\mu_{ts} - (\mu_{tu}+\mu_{us})\big| \leq c_1 |t-s|^\zeta,
$$
for some positive constant $c_1$ and exponent $\zeta>1$, then there exists a unique path $x_\bullet$ started from $0$ whose increments $x_{ts}$ satisfy
$$
\big|x_{ts}-\mu_{ts}\Big| \lesssim |t-s|^\zeta.
$$
This is Gubinelli-Feyel de la Pradelle' sewing lemma \cite{Gubinelli, FdlP}.

\ssk

Now, assume we are given an $\RR^\ell$-valued (weak geometric) $\alpha$-H\"older rough path 
$$
{\bfX} = \big((X_{ts},\bbX_{ts})\big)_{0\leq s\leq t\leq T},
$$ 
with $X_{ts}\in\RR^\ell$ and $\bbX_{ts}\in \RR^\ell\otimes\RR^\ell$, and a map $\textrm{F}\in C^2\Big(\RR^d,\textrm{L}\big(\RR^\ell,\RR^d\big)\Big)$. Following Lyons, an $\RR^d$-valued path $x_\bullet$ is said to solve the rough differential equation 
\begin{equation}
\label{EqRDE}
dx_t = \textrm{F}(x_t)\,{\bfX}(dt)
\end{equation}
if one has 
\begin{equation}
\label{EqRDEdefinition}
x_t-x_s = \textrm{F}(x_s)\,X_{ts} + \textrm{F}'(x_s)\textrm{F}(x_s)\,{\bbX}_{ts} + O\big(|t-s|^\zeta\big)
\end{equation}
for all $0\leq s\leq t\leq T$, for some constant $\zeta>1$. Gubinelli's crucial remark in \cite{Gubinelli} was to notice that for a path $x_\bullet$ to satisfy equation \eqref{EqRDEdefinition}, it needs to be controlled by (the first level) $X$ (of the rough path $\bfX$) in the sense that one has
\begin{equation*}
\label{EqTaylorOneRDE}
x_t-x_s = x'_s\,X_{ts} + O\big(|t-s|^{2\alpha}\big),
\end{equation*}
for some $\textrm{L}(\RR^\ell,\RR^d)$-valued $\alpha$-H\"older path $x'_\bullet$, here $x'_s = \textrm{F}(x_s)$. 

\begin{definition}
A $\textrm{\emph{V}}$-valued \textbf{path} $(z_t)_{0\leq t\leq T}$ is said to be \textbf{controlled by} $\bfX$ if its increments $Z_{ts} = z_t-z_s$, satisfy 
$$
Z_{ts} = Z'_sX_{ts} + R_{ts},
$$
for all $0\leq s\leq t\leq T$, for some $\textrm{\emph{L}}(\textrm{\emph{E}},\textrm{\emph{V}})$-valued $\frac{1}{p}$-Lipschitz map $Z'_\bullet$, and some $\textrm{\emph{V}}$-valued $\frac{2}{p}$-Lipschitz map $R$. The controlled path itself in not the path $z_\bullet$ but rather the pair $(Z,Z')$. One endows the set $\mcC_\bfX(\textrm{\emph{V}})$ of such paths with a Banach space structure setting 
$$ 
\big\|(Z,Z')\big\| := \big\|Z'\big\|_\frac{1}{p} + \|R\|_\frac{2}{p} + |z_0|.
$$
\end{definition}

The point of this notion is that, somewhat conversely, if we are given an $\textrm{L}(\RR^\ell,\RR^d)$-valued $\alpha$-H\"older path $a_\bullet$ controlled by $\bfX$, then there exists, by the sewing lemma, a unique $\RR^d$-valued path $b_\bullet$ whose increments satisfy 
$$
b_t-b_s = a_s\,X_{ts} + a'_s\,\bbX_{ts} + O\big(|t-s|^\zeta\big),
$$
for some exponent $\zeta>1$, because the formula
$$
\mu_{ts} := a_s\,X_{ts} + a'_s\,\bbX_{ts}
$$
defines an almost additive $2$-index map. With a little bit of abuse, we write $\int_0^\bullet a_s\,{\bfX}(ds)$ for that path $b_\bullet$ -- this path depends not only on $a$ but rather on  $(a,a')$. This path depends continuously on $(a,a')\in\mcC_{\bfX}(\textrm{V})$ and $\bfX$. 

\medskip

Now, if we are given a rough path $\bfX$ over E, defined on the time interval $[0,1]$, a one form $\textrm{F}\in C^k\Big(\textrm{V},\textrm{L}\big(\textrm{U},\textrm{V}\big)\Big)$, with $k\geq 3$, and controlled paths $z\in\mcC_\bfX(\textrm{V})$ and $y\in\mcC_\bfX(\textrm{U})$, one defines the rough integral 
$$ 
\int_0^\bullet \textrm{F}(z_u)\,dy_u
$$
as the unique additive functional uniquely associated with the almost-additive map
$$
\textrm{F}(z_s)y_{ts} + \big(D_{z_s}\textrm{F}\big)(Z'_s\otimes y'_s)\,\bbX_{ts}. 
$$

\medskip

Denote by $\frak{I}(\textrm{F},y)$ the It\^o map that associates to any V-valued $\textrm{Lip}_k$ one form F on U, and any path $y_\bullet$ in U controlled by $\bfX$, the unique solution of the rough differential equation 
\begin{equation}
\label{EqRDE}
x_\bullet = x_0 + \int_0^\bullet \textrm{F}(x_r)\,dy_r, 
\end{equation}
with given initial condition $x_0\in \textrm{V}$. (See for instance \cite{Lyons97,LyonsStFlour,BailleulInfiniteDim, FH13, Gubinelli} for references on rough differential equations in an infinite dimensional setting. Solving such a rough differential equation is simply finding a fixed point for the above well-defined integral equation.) The It\^o map is defined on the product Banach space  
$$
(\textrm{F},y_\bullet)\in\textrm{L}\big(\textrm{U},\textrm{Lip}_k(\textrm{V},\textrm{V})\big)\times\mcC_{\bfX}(\textrm{U})
$$ and takes values in the affine subspace $\mcC^{x_0}_{\bfX}(\textrm{V})$ of $\mcC_{\bfX}(\textrm{V})$, made up of V-valued paths controlled by $\bfX$ started from $x_0$; so there is no difficulty in defining Fr\'echet differentiability in that setting. Despite the simplicity of its proof, the following regularity result is our main result.

\begin{thm}
\label{ThmRegularity}
Fix $x_0\in \textrm{V}$. The It\^o-Lyons solution map 
$$
\frak{I} : \textrm{\emph{L}}\big(\textrm{\emph{U}},\textrm{\emph{Lip}}_k(\textrm{\emph{V}},\textrm{\emph{V}})\big)\times\mcC_{\bfX}(\textrm{\emph{U}})\rightarrow \mcC^{x_0}_{\bfX}(\textrm{\emph{V}})
$$ 
is $\mcC^{[k]-2}$ Fr\'echet-differentiable.
\end{thm}

\medskip 

A similar regularity result has been proved by Li and Lyons in \cite{LiLyons} for the solution map to a Young differential equation with fixed F, corresponding to the case $1<p<2$. Their result is more general than the statement corresponding in that setting to theorem \ref{ThmRegularity}, for the It\^o map they consider is defined on a much bigger space than $\mcC_\bfX(\textrm{U})$; this is the reason why the proof of their result requires non-trivial arguments. On the other hand, the non-linear nature of the space of rough paths, for $2\leq p<3$, makes the extension of their approach to that setting hard, as illustrated by the work \cite{QianTudor} of Qian and Tudor, in which they introduce an appropriate notion of tangent space to the space of $p$-rough paths which is of delicate use. The present use of the Banach setting of controlled paths somehow linearizes many considerations, while providing results that are still widely applicable. Note that one cannot hope for boundedness of the derivative of $\frak{I}$ in the above level of generality; so the It\^o-Lyons solution map is a priori only locally Lipschitz.

\bigskip

\begin{DemThmRegularity}
Denote by $\mcC^0_{\bfX}(\textrm{V})$ the space of V-valued paths controlled by $\bfX$ started from $0$. Our starting point is the fact that the solution to equation \eqref{EqRDE} is the unique fixed point to the map
$$
\Phi(\textrm{F},y,x_0) : (z,z') \mapsto \left(\int_0^\bullet \textrm{F}(x_0 + z_r)\,dy_r, \textrm{F}(x_0 + z_\bullet)\,y'_\bullet\right) 
$$
from $\mcC^0_{\bfX}(\textrm{V})$ to itself. Write $\mcC^0_{{\bfX},T}(\textrm{V})$ for the restriction to the time interval $[0,T]$ of the controlled paths in $\mcC^0_{\bfX}(\textrm{V})$. As shown in proposition 8 of Gubinelli's work \cite{Gubinelli} there exists a time $T\leq 1$, independent of $x_0$, and uniformly bounded away from $0$, for F and $y$ in bounded sets of their respective spaces, such that $\Phi(\textrm{F},y,x_0)$ is a contraction of a fixed ball in $\mcC^0_{{\bfX},T}(\textrm{V})$ that contains all possible solutions to equation \eqref{EqRDE}, by a priori estimate such as those proved in section 5 of \cite{Gubinelli}.
 
\ssk

Note that the map $\Phi$ depends linearly and continuously on  $y$ and F by classical estimates on rough integrals, as can be found for instance in theorem 4.9 in \cite{FH13}, so it is in particular a smooth function of $y_\bullet\in\mcC_{\bfX}(\textrm{U})$. It is also elementary to see that the map $\Phi$ is a $\mcC^{[k]-2}$ function of all its arguments -- the dependence on $x_0$ can be treated technically as a dependence on F. Write $\partial_{(z,z')}\Phi$ for the derivative of $\Phi$ with respect to $(z,z')$. Since $\Phi(\textrm{F},y,x_0)$ is uniformly strictly contracting on a fixed ball of $\mcC^0_{{\bfX},T}(\textrm{V})$, for F and $y$ ranging in some bounded sets, independently of $x_0\in\RR^d$, the continuous linear map $\textrm{Id}-\partial_{(z,z')}\Phi(\textrm{F},y,x_0\, ;\,\cdot)$ from $\mcC^0_{{\bfX},T}(\textrm{V})$ to itself has an inverse given under the classical form of the converging Neumann series 
$$
\sum_n \big\{\partial_{(z,z')}\Phi(\textrm{F},y,x_0\,;\,\cdot)\big\}^n.
$$
The map $\textrm{Id}-\partial_{(z,z')}\Phi(\textrm{F},y,x_0\, ;\,\cdot)$ is thus a continuous isomorphism of $\mcC^0_{{\bfX},T}(\textrm{V})$, by the \textit{open mapping theorem}. It is then a direct consequence of the \textit{implicit function theorem} that the unique fixed point $x_\bullet$ to the equation 
$$
\Phi(\textrm{F},y,x_0\, ;\,z)=z
$$
in $\mcC^0_{{\bfX},T}(\textrm{V})$ defines a $\mcC^{[k]-2}$ function of F and $y$ and $x_0$. An elementary patching procedure can be used to extend the result to the whole space $\mcC^0_{\bfX}(\textrm{V})$. 
\end{DemThmRegularity}

\bigskip

It should be clear to the reader that the above scheme of proof extends in a straightforward way to deal with rough differential equations driven by more irregular paths. Gubinelli's branched rough paths \cite{GubinelliBranched} provide a convenient setting for that purpose; we leave the details to the reader.  As another remark, note that if $\bfX$ is a piece of a "higher dimensional" rough path $\widehat\bfX$, one can work in $\mcC_{\widehat{\bfX}}(\textrm{U})\supset \mcC_\bfX(\textrm{U})$  and perturb the original rough differential equation \eqref{EqRDE} with rough signals non-controlled by $\bfX$, but rather controlled by $\widehat{\bfX}$, that is controlled by other "pieces" of $\widehat{\bfX}$. More concretely, let $\widehat{\bfX}$ stand for instance for the canonical rough path above a $(d_1+d_2)$-dimensional fractional Brownian motion with Hurst index $\frac{1}{3}< H\leq \frac{1}{2}$, and let $\pi_i$ stand for the canonical projection from $\RR^{d_1+d_2}$ to $\RR^{d_i}$, for $i\in\{1,2\}$, extended canonically to the $2$-level truncated tensor space over $\RR^{d_1+d_2}$. Then ${\bfX} := \pi_1\big(\widehat{\bfX}\big)$ is the canonical rough path above the $d_1$-dimensional fractional Brownian motion. Write here $\widehat{\frak{I}}$ for the It\^o map 
$$
\widehat{\frak{I}} : \textrm{L}\big(\textrm{U},\textrm{Lip}_k(\textrm{V},\textrm{V})\big)\times\mcC_{\widehat \bfX}(\textrm{U})\rightarrow \mcC^{x_0}_{\widehat\bfX}(\textrm{V}). 
$$
Since we have
$$
\widehat{\frak{I}}(\textrm{F},y_\bullet) = \frak{I}(\textrm{F},y_\bullet),
$$
for $y_\bullet\in\mcC_{\bfX}(\textrm{U})$,  our framework allows us to differentiate the solution to the rough differential equation \eqref{EqRDE} in the direction of $\pi_2\big(\widehat{\bfX}\big)$, as is useful in the setting of Malliavin calculus, even in the rough paths setting -- see the work \cite{InahamaMalliavin} of Inahama for an illustration of this point. See also section \ref{SubSectionTaylorExpansion} for a related point.

\bigskip

\section{Applications}
\label{SectionApplications}

We illustrate in this section the use of theorem \ref{ThmRegularity} with two typical examples: Taylor-like expansions for solutions of rough differential equations, as studied for instance in Aida \cite{Aida}, Inahama-Kawabi \cite{InahamaKawabi} and Inahama \cite{Inahama}, and the construction of some natural dynamics on the path path space over a compact Riemannian manifold, such as Driver's flow. We also describe in the last section how the above regularity theorem can be used to prove an integration by parts formula for the It\^o-Lyons map.

\medskip

\subsection{Taylor expansion of Inahama-Kawabi type}
\label{SubSectionTaylorExpansion}

Inahama and Kawabi developed in \cite{InahamaKawabi} and \cite{Inahama} some sophisticated machinery to obtain some Taylor-like expansion in $\ep$ for solutions of rough differential equations of the form 
 \begin{equation}
\label{EqZep} 
dz^{(\ep)}_t = \sigma\big(\ep,z^{(\ep)}_t\big) d{\bfX}_t + b\big(\ep,z^{(\ep)}_t\big)d\Lambda_t, 
 \end{equation}
 where $\bfX$ is a given infinite dimensional $p$-rough path over U, and $\Lambda$ a path with finite $q$-variation, with $\frac{1}{p}+\frac{1}{q}>1$, and with values in some other Banach space $\textrm{U}'$. These estimates were used in an instrumental way to obtain some asymptotic formulas for quantities of the form
$$ 
\EE\left[A\big(z^{(\ep)}_\bullet\big)\exp^{-\frac{1}{\ep^2}B(z^{(\ep)}_\bullet)}\right], 
$$ 
for some real-valued function $A,B$ on $\textrm{U}$, using Laplace method, as extended first to a finite dimensional stochastic setting by Azencott \cite{Azencott} and Ben Arous \cite{BenArous}. Such Taylor expansions 
$$
Z^{(\ep)}_\bullet = Z^0_\bullet + \ep Z^1_\bullet + \cdots + \ep^{[k]-2} Z^{[k]-2}_\bullet  + o\Big(\ep^{[k]-2}\Big)
$$
are direct consequences of theorem \ref{ThmRegularity}, in the driftless case where $\Lambda = 0$, provided the map
$$
\ep \mapsto \sigma(\ep,\cdot) \in\textrm{L}\big(\textrm{U},\mcC^k(\textrm{V},\textrm{V})\big)
$$ 
is at least $\mcC^{[k]-2}$. The path $Z^0$ is null in that case and a formal expansion of equation \wrt $\ep$ provides the set of time-dependent affine equations satisfied by the $Z^i$'s; see section 4 of \cite{Inahama} for instance. One can give a similar expansion when the drift is non-null by constructing first the canonical rough path $\widehat\bfX$ over $\textrm{U}\times\textrm{U}'$ above $\bfX$ and $\Lambda$ (recall $\frac{1}{p}+\frac{1}{q}>1$), and then work in $\mcC_{\widehat\bfX}(\textrm{V})$ where theorem \ref{ThmRegularity} applies. This approach to Taylor expansion for solutions of equation \eqref{EqZep} somehow simplifies the orginal approach of Inahama and Kawabi.

\medskip

\subsection{Dynamics on path space}
\label{SubSectionDynamicsPathSpace}

Assume here U=V is a given Banach space and F is a fixed $\textrm{Lip}_3$-one form on U; fix an initial condition $x_0\in\textrm{U}$. As the It\^o-Lyons map 
$$
\frak{I}(\textrm{F},\cdot) : \mcC^{x_0}_{\bfX}(\textrm{U}) \rightarrow\mcC^{x_0}_{\bfX}(\textrm{U})
$$ 
is locally Lipschitz, one can apply the classical Cauchy-Lipschitz theorem in $\mcC^{x_0}_{\bfX}(\textrm{V})$ and solve uniquely the ordinary differential equation
$$
\frac{dy_\bullet(t)}{dt} = \frak{I}\big(\textrm{F}, y_\bullet(t)\big)
$$
for a given initial condition $y_\bullet$ in $\mcC^{x_0}_{\bfX}(\textrm{V})$, on a maximal interval of definition. This simple remark provides a straightforward approach to the main results of \cite{LyonsQianFlow} and can be used in a number of geometrical situations to define some natural evolutions on the path space over some given manifold. (Note that this evolution on path space keeps the starting point fixed.) We give two representative examples of such dynamics below, of which Driver's flow \cite{Driver} on the path space over a compact manifold is an instance. 

\medskip

\subsubsection{Driver's flow equation}

Let $\MMM$ be a compact smooth $n$-dimensional submanifold of $\RR^d$, endowed with the Riemannian structure inherited from the ambiant space. Nash's theorem ensures there is no loss of generality in considering that setting. Let $\nabla : \RR^d\times\mcC^\infty\big(\RR^d,\RR^d\big)\rightarrow\mcC^\infty\big(\RR^d,\RR^d\big)$, be a compactly supported smooth extension to $\RR^d$ of the Levi-Civita covariant differentiation operator on $T\MMM$. This $\RR^d$-dependent operator has a natural extension to sections of $\textrm{L}(\RR^d)$, still denoted by the same symbol. 

\ssk

Denote by $\mcC_\bfX(\MMM)$ the subset of $\mcC_\bfX(\RR^d)$ made up of $\MMM$-valued paths $y_\bullet$ controlled by $\bfX$. Given $T_0\in\textrm{L}(\RR^d)$, the solution to the rough differential equation in $\textrm{L}(\RR^d)$
\begin{equation}
\label{EqParallelTransport}
dT_s = \nabla(T_s)\,dy_s,
\end{equation}
started from $T_0$, defines a path in $\textrm{L}(\RR^d)$ above $y_\bullet$. If the path $y_\bullet$ takes values in $\MMM$ and $T_0$ is orthonormal and has its first $n$ columns forming an orthonormal basis of $T_{y_0}\MMM$, then the first $n$ columns of $T_s$ form an orthonormal basis of $T_{y_s}\MMM$ at any time $0\leq s\leq 1$. Also, the restriction of $\big(T_s\big)_{0\leq s\leq 1}$ to $\RR^n$ does not depend on ${T_0}_{\big|(\RR^n)^\perp}$, and $f_h$ depends only on the restriction to $\MMM$ of $\nabla$; see for instance the very recent work \cite{CassDriverLitterer} of Cass, Driver and Litterer on constrained rough paths, section 5. The map $T_s$ transports parallelly $T_{y_0}\MMM$ to $T_{y_s}\MMM$ as an isometry. In that case, and given any controlled path $h_\bullet$ in $\mathcal{C}_{\bfX}(\RR^n)$, with $\RR^n$ seen as a subset of $\RR^d$, we define a Lipschitz continuous map from $\mcC_{\bfX}(\MMM)$ to $\mcC_{\bfX}(\RR^d)$ setting
$$
\frak{F}_h(y)_s = T_sh_s, 
$$
for $0\leq s\leq 1$. The vector $f_h(y_\bullet)_s$ belongs to $T_{y_s}\MMM$ at any time $0\leq s\leq 1$. It is elementary to see that if $y_\bullet\in\mcC_{\bfX}(\MMM)$, then the solution in $\mcC_{\bfX}(\RR^d)$ to the ordinary differential equation
\begin{equation}
\label{EqPathSpaceEvolution}
\frac{dy_\bullet(t)}{dt} = \frak{F}_h\big(y_\bullet(t)\big)
\end{equation}
started from $y_\bullet$ takes values in $\mcC_\bfX(\MMM)$. The local flow associated with the above dynamics is precisely Driver's flow when $h\in \textrm{H}^1(\RR^n)$ and $\bfX$ is a typical realization of the $n$-dimensional Brownian rough path. As Driver's flow is known to be almost-surely defined for all times, the local flow defined by equation \eqref{EqPathSpaceEvolution} is actually almost-surely globally defined. This construction of Driver's flow simplifies the corresponding construction by Lyons and Qian given in \cite{LyonsQianFlow}, where a local flow on the space of $\MMM$-valued paths with finite $1$-variation is constructed and then extended to a rough path space by a continuity argument. The present approach illustrates the interest of working in the Banach setting of controlled paths, as opposed to the non-linear setting of rough paths.

\ssk

Note that we do not really need a Riemannian metric on $\MMM$ for the only object we use is a connection on $T\MMM$. Of course, the main interest of Diver's flow in a Riemannian/Brownian setting, and the main reason for its introduction, is the fact that it leaves Wiener measure on the path space over $\MMM$ quasi-invariant, provided $h$ satisfies some conditions; see \cite{Driver} and \cite{Hsu}. We do not touch upon this point here as the theory of rough paths does not seems to bring any insight on that question.

\bigskip

\subsubsection{Varying the connection}

As another example of use of theorem \ref{ThmRegularity} for constructing natural dynamics on the path space over a compact Riemannian manifold, one can follow Lyons and Qian \cite{LyonsQianJacobi} and look at the flow generated by the first order variation of a parameter-dependent connection on $T\MMM$. Its construction goes as follows; as above, we consider $\MMM$ as a submanifold of $\RR^d$, whose Riemannian structure is inherited from the ambiant flat metric on $\RR^d$. The reference rough path $\bfX$ is fixed.

\medskip

Let P be a smooth map on $\RR^d$, with bounded derivatives and values in the set of projectors, such that $\textrm{P}(x)$ has range $T_x\MMM$ at any point $x$ of $\MMM$. Setting 
$$
\big(\nabla^\textrm{P}_{\sf u} {\sf v}\big)(x) = \textrm{P}(x)\Big\{\big(D_x {\sf v}\big)\big({\sf u}(x)\big)\Big\},
$$ 
for any two vector fields ${\sf u,v}$ on the ambiant space $\RR^d$, we get back for instance the Levi-Civitta connection on $\MMM$ if $\textrm{P}(x)$ is the orthogonal projection on $T_x\MMM\subset\RR^d$ at any point $x$ of $\MMM$ and the vector fields ${\sf u}$ and ${\sf v}$ are tangent to $\MMM$ on $\MMM$. Given any other such map Q, and some $\MMM$-valued controlled path $y_\bullet\in\mcC_{\bfX}(\MMM)$, let first $z_\bullet$ be the anti-development in $T_{y_0}\MMM$ of $y_\bullet$. It is defined as the solution of the rough differential equation
$$
dz_s = T_s^{-1}\textrm{P}(y_s)\,dy_s
$$
started from $0\in T_{y_0}\MMM$, where $T_s$ is the parallel transport map constructed in equation \eqref{EqParallelTransport}; so the pair $(T_\bullet,z_\bullet)$ is solution to a rough differential equation driven by $y$ and smooth 1-forms, so it is a smooth function of $y\in\mcC_{\bfX}(\RR^d)$. We refer the reader to the above mentionned work \cite{CassDriverLitterer} for more explanations on this anti-development map. Let now $w^\ep_\bullet$ be the solution to the rough differential equation
$$
dw^\ep_s = (\textrm{P}+\ep\textrm{Q})\big(w^\ep_s\big)\,dz_t
$$
started from $y_0\in\MMM$, so $w^\ep = \frak{I}(\textrm{P}+\ep\textrm{Q},z)$. The formula
$$
\frak{F}(y) = \frac{d}{d\ep}_{\big| \ep=0} \,\frak{I}(\textrm{P}+\ep\textrm{Q},z)
$$
defines a controlled path which takes values in $T_{y_s}\MMM$ at any time $0\leq s\leq 1$; as a function of $y_\bullet$, it is locally Lipschitz continuous, so it generates a local flow on the path space over $\MMM$.

\ssk

Note that this construction works regardless of any assumption on the connectors P and Q, which goes farther than the metric and Riemannian setting adopted in \cite{LyonsQianFlow}.

\bigskip

\subsection{Integration by parts for the It\^o-Lyons map}
\label{SubsectionIBP}
 
As a last illustration of the use of the regularity theorem, theorem \ref{ThmRegularity}, we explain in this section how to prove an integration by parts formula for functionals of some particular classes of random rough paths, leaving the details to the reader. 

\ssk 

We restrict our attention in this paragraph to the case where $V=\RR^\ell$ and $U=\RR^d$ have finite dimension, so $\textrm{F}=(V_1,\dots,V_\ell)$ is a collection of $\ell$ vector fields on $\RR^d$, of class $\mcC^5_b$, and consider the equation
\begin{equation}
\label{EqRDEFinite}
dx_t = \textrm{F}(x_t){\bfX}(dt). 
\end{equation}
Denote by 
$$
\Gamma = \int_0^1 \big(U_r^{-1}V_i\big)(x_r)\otimes \big(U_r^{-1}V_i\big)(x_r)dr
$$
the usual Malliavin covariance symmetric matrix, where $(U_r)_{0\leq r\leq 1}$ stands for the derivative flow of equation \eqref{EqRDE}, with values in $\textrm{GL}\big(\RR^d\big)$. Identify $\Gamma$ with its associated symmetric bilinear form. The following integration by parts result is proved below as a corollary of theorem \ref{ThmRegularity}.

\begin{thm}[Integration by parts formula]
\label{ThmIBP}
Let  $2\leq p<3$, and $\bfX$ be the rough path above the $\RR^\ell$-valued fractional Brownian motion with Hurst index $\frac{1}{p}$. Then there exists an $\bfX$-dependent second order differential operator $L^\bfX$ such that we have 
$$
\EE\Big[\Gamma\big(\nabla f(x_1),\nabla g(x_1)\big)\Big] = -2\,\EE\Big[f(x_1)\,\big(L^{\bfX}g\big)(x_1)\Big],
$$
for any real-valued functions $f,g$ on $\RR^d$ of class $\mcC^3_b$.
\end{thm}

As is well-known, this formula is the cornerstone to deriving absolute continuity results for the distribution of $x_1$. So far, all the works on this subject, by Cass, Friz \cite{CassFriz} , Hairer and co \cite{CassHairerLittererTindel}, used a result of Bouleau-Hirsch to link the regularity properties of the distribution of $x_1$ to the almost-sure and integrability properties of the Malliavin covariance matrix. Theorem \ref{ThmIBP} provides a pure rough path approach to the subject, in the above setting. The sketch of proof given below makes it clear that the above integration by parts formula holds for a wider class of random rough paths, whose laws are invariant by the action of the isometries of space, and for which one has strong integrability  estimates such as those proved by Cass-Hairer-Litterer-Tindel in \cite{CassHairerLittererTindel}. 

\bigskip 

The idea of the proof is simple and goes back to Malliavin's original approach \cite{Malliavin}. (See also the work \cite{NorrisTwistedSheetMalliavin} of Norris for a related approach using two-parameter semimartingales.) As the one form F will be fixed here, write $\frak{I}_1(\bfX)$ for the time $1$ value of the solution to equation \eqref{EqRDEFinite}, with a slight abuse of notations. Start by recalling that $\textrm{L}\big(\RR^\ell\big)$ has a natural action on the set of $p$-rough paths; write $\textrm{M}\bfX$ for the action of $\textrm{M}\in\textrm{L}\big(\RR^\ell\big)$ on $\bfX$. Assume that the distribution of $\bfX$ is invariant by the action of the group $\textrm{O}\big(\RR^\ell\big)$ of isometries of $\RR^\ell$. Let $(W_s)_{0\leq s\leq 1}$ be a $1$-dimensional Brownian motion with initial position uniformly distributed on $[0,2\pi)$, and $\textrm{M}^s\in \textrm{O}\big(\RR^\ell\big)$ be the rotation of angle $W_s$ in the $2$-plane generated by the first two vectors of the canonical basis. Set 
$$
{\bfX}^s = \textrm{M}^s\bfX.
$$
The solution $x^s_\bullet$ to the equation 
$$
dx^s_t = \textrm{F}(x_t){\bfX}^s(dt)
$$
is controlled by $\bfX$, for every $0\leq s\leq 1$, and the basic reversibility property 
$$
\big({\bfX}^0,{\bfX}^s\big) \overset{\textrm{law}}{=} \big({\bfX}^s,{\bfX}^0\big)
$$
holds for all $0\leq s\leq 1$. It follows from the latter that 
{\small
\begin{equation*}
\begin{split}
\EE\Big[\Big\{f\big(\frak{I}_1({\bfX}^s)\big) - f\big(\frak{I}_1({\bfX}^0)\big)\Big\} &\Big\{g\big(\frak{I}_1({\bfX}^s)\big)-g\big(\frak{I}_1({\bfX}^0)\big)\Big\}\Big] \\ 
&= -2 \EE\Big[f\big(\frak{I}_1({\bfX}^0)\big)\,\Big\{g\big(\frak{I}_1({\bfX}^s)\big)-g\big(\frak{I}_1({\bfX}^0)\big)\Big\}\Big]
\end{split}
\end{equation*}
}
As we are working in the fixed space $\mcC_{\bfX}(\RR^d)$ along the $s$-deformation, one can divide both sides by $s$ and make a Taylor expansion using theorem \ref{ThmRegularity}. Writing $s=\sqrt{s}\sqrt{s}$ on the left hand side, we capture the linear part of $f\big(\frak{I}_1({\bfX}^s)\big)-f\big(\frak{I}_1({\bfX}^0)\big)$ \wrt $W_s$, with negligeable remainders. The $\Gamma$ term comes from 
the derivative flow of equation \eqref{EqRDE}. The terms that are relevant on the right hand side come from a second order expansion of $g\big(\frak{I}_1({\bfX}^s)\big)$, from which the above mentionned $\bfX$-dependent second order differential operator comes from. Integrability estimates like those proved in \cite{CassHairerLittererTindel} by Cass-Hairer-Litterer-Tindel are needed to justify the use of Lebesgue's dominated convergence while sending $s$ to $0$. We leave the reader the pleasure to fill in the details.

\vfill

\end{document}